\documentclass[12pt]{amsart}
\usepackage{amssymb,amsmath}
\makeatletter
\def\@cite#1#2{{\m@th\upshape\bfseries%
[{#1\if@tempswa{\m@th\upshape\mdseries, #2}\fi}]}}
\makeatother
\theoremstyle{plain}
\newtheorem{thm}{Theorem}[section]
\newtheorem{lem}[thm]{Lemma}

\newtheorem{prop}[thm]{Proposition}
\theoremstyle{definition}
\newtheorem{rem}[thm]{Remark}
\newtheorem{defn}[thm]{Definition}
\newtheorem{note}[thm]{Note}
\newtheorem{notes}[thm]{Notes}
\newtheorem{eg}[thm]{Example}

\newtheorem{prob}[thm]{Problem}
\newcommand{\Prf}{\noindent\textbf{Proof.\ }}
\newcommand{\bx}{\strut\hfill$\blacksquare$\medbreak}

\newcommand{\ca}{\mathrm{C}^*}

\newcommand{\ol}{\overline}

\DeclareMathOperator*{\sotlim}{\textsc{sot}--lim}

\newcommand{\wot}{\textsc{wot}}

\newcommand{\bbA}{{\mathbb{A}}}

\newcommand{\bbC}{{\mathbb{C}}}

\newcommand{\bbF}{{\mathbb{F}}}
\newcommand{\bbG}{{\mathbb{G}}}

\newcommand{\bbT}{{\mathbb{T}}}

\newcommand{\bbZ}{{\mathbb{Z}}}
 \newcommand{\B}{{\mathcal{B}}}

\renewcommand{\H}{{\mathcal{H}}}
 
 \newcommand{\J}{{\mathcal{J}}}
 \newcommand{\K}{{\mathcal{K}}}  

 \newcommand{\M}{{\mathcal{M}}}
 
\renewcommand{\O}{{\mathcal{O}}}
\renewcommand{\P}{{\mathcal{P}}}

\newcommand{\upchi}{{\raise.35ex\hbox{$\chi$}}}
\newcommand{\fA}{{\mathfrak{A}}}

\newcommand{\fS}{{\mathfrak{S}}}

\newcommand{\qand}{\quad\text{and}\quad}

\newcommand{\qfor}{\quad\text{for}\quad}

\newcommand{\qiff}{\quad\text{if and only if}\quad}
\newcommand{\Alg}{\operatorname{Alg}}

\newcommand{\spn}{\operatorname{span}}

\newcommand{\BH}{\B(\H)}
\newcommand{\fnplus}{\bbF_n^+}
\newcommand{\fom}{\bbF_\omega}
\newcommand{\sumin}{\sum_{i=1}^n}

\newcommand{\rowt}{(T_1, \ldots, T_n)}
\newcommand{\rows}{(S_1, \ldots, S_n)}

\newcommand{\fu}{{\bbF}_\omega}
\newcommand{\fock}{\ell^2 ({\bbF}_\omega)}

\newcommand{\fsulam}{\fS_{\omega,\Lambda}}
\newcommand{\munegk}{\mu_{\omega_m^{-1}}}
\newcommand{\qri}{Q_{r,i}}
\newcommand{\pri}{P_{r,i}}

\newcommand{\ju}{\J_{u,k}}

\newcommand{\uomk}{U_{\omega,k}}
\newcommand{\fomk}{\bbF_{\omega,k}}

\newcommand{\fomkl}{\bbF_{\omega,k}^{(l)}}
\newcommand{\omprojns}{\omega_m(S)\omega_m(S)^*}

\begin{document}

\title[Noncommutative Bilateral Weighted Shifts
]%
{On bilateral weighted shifts in noncommutative multivariable operator 
theory}
%
\author[D.W.Kribs]{David~W.~Kribs${}^1$} 
\address{Department of Mathematics and Statistics, University of Guelph, Guelph, ON,
CANADA \\ N1G 2W1}
\email{kribs@math.purdue.edu}
\thanks{{\it key words and phrases. } Hilbert space, operator, bilateral 
weighted 
shift, periodicity, reducing subspaces, infinite word, noncommutative 
multivariable operator theory,  
nonselfadjoint operator algebras, Fock space. }
\thanks{${}^{1\;}$partially supported by a Canadian NSERC Post-doctoral
Fellowship.}

\date{}
\begin{abstract}
We present a generalization of bilateral weighted shift operators for the 
noncommutative multivariable setting. 
We discover a notion of periodicity for these shifts, which has an 
appealing diagramatic interpretation in terms of an infinite tree 
structure associated with the underlying Hilbert space. 
These shifts 
arise naturally through weighted versions of certain representations  of 
the Cuntz $\ca$-algebras $\O_n$. It is convenient, and equivalent, to 
consider the weak operator topology closed algebras generated by these 
operators when investigating their joint reducing subspace structure. 
We prove these 
algebras have non-trivial reducing subspaces exactly when the shifts are 
doubly-periodic; that is, the weights for the shift have periodic 
behaviour, and the corresponding representation of $\O_n$ has a certain 
spatial periodicity.  This generalizes Nikolskii's Theorem 
for the single 
variable case. 
\end{abstract}

\maketitle

In \cite{nonsa} and \cite{inductive}, we began studying versions of 
unilateral weighted shift operators in noncommutative multivariable 
operator theory. We called them {\it weighted shifts on Fock space} since 
they act naturally on the full Fock space Hilbert space. 
These shifts and the algebras they generate were first studied by Arias and Popescu \cite{AP2}
from the perspective of {\it weighted Fock spaces}. 
The basic goals  of this program 
are to extend results from the commutative (single variable) setting and, 
at the same time,  
expose new noncommutative phenomena. 
In the current paper, we continue this line of investigation by presenting 
versions of {\it bilateral} weighted shift operators for the 
noncommutative multivariable setting. Our analysis is chiefly spatial in 
nature: we examine the joint reducing 
subspace structure for these shifts, and consider 
reducibility questions for the weak operator topology closed 
(nonselfadjoint) algebras they 
generate. In particular, we give a complete characterization of 
reducibility 
for the 
algebras strictly in terms of two notions of periodicity 
associated 
with these shifts. This generalizes Nikolskii's Theorem 
\cite{Nikol} for reducing subspaces of bilateral weighted shift operators 
on Hilbert space. 

In the first section we include a short introduction to the subject. Next 
we present the definition and derive some basic properties for these 
shifts.  The second section contains the main results of 
the paper; most importantly, a complete description of the reducing 
subspace structure for the shifts and the operator algebras they 
generate. In the final section we include examples and state some open problems.

\section{Introduction}\label{S:intro}

Given a positive integer $n\geq 1$, or $n= \infty$, we may consider all 
possible $n$-tuples of operators 
$S= \rows$ which act on a common Hilbert space $\H$ and satisfy the 
relations 
\begin{eqnarray}\label{cuntz}
S_i^*S_i = I \qfor 1\leq i \leq n \qand \sumin S_iS_i^* = 
I.
\end{eqnarray}
When $n=1$, the single variable case, we are simply talking about unitary 
operators. Distinguished amongst the collection of unitary operators on  
(separable) infinite-dimensional spaces is 
the canonical {\it 
bilateral shift} which acts on an 
orthonormal basis $\{ \xi_k : k\in\bbZ\}$ for $\H$ by $U\xi_k = 
\xi_{k+1}$. 
Bilateral shifts play a central role in operator theory generally, and 
specifically within the structure theory for isometries. Thus, bilateral 
{\it weighted} shifts, operators $T$ defined on $\H$ by $T\xi_k = 
\lambda_{k+1}\xi_{k+1}$, $\lambda_{k+1}\in\bbC$, have also been studied 
for some time, often providing new classes of examples and shedding light 
on the 
theory overall (see the survey article \cite{Shields} for example). 

On the other hand, the set of all $n$-tuples $S= \rows$ which satisfy 
(\ref{cuntz}) for $n\geq 2$ is extremely vast, and there perhaps may not 
be an analogous structure theory here. For instance, the universal 
$\ca$-algebra generated by (\ref{cuntz}) (that is, the norm closure of the 
orthogonal direct sum of all possible $\ast$-representations determined by 
$n$-tuples satisfying (\ref{cuntz})), called the Cuntz algebra $\O_n$ 
\cite{Cun}, has such a rich representation theory that it is not possible 
to classify all its representations up to unitary equivalence. (It is an 
NGCR algebra \cite{Glimm}.) 
Thus, there has been considerable recent interest in studying various 
subclasses of representations of $\O_n$ for a number of reasons. For 
examples from some of the diverse perspectives we 
reference, in passing,  work of Ball, Bratteli, Davidson, Jorgensen, 
Katsoulis, Pitts, Shpigel, Vinnikov, and the author 
\cite{BV1,BV2,BJit,BJwave,Dsurvey,DKP,DKS,DP1,quantum}. 

While there is not a well developed structure theory for $n\geq 2$, there 
are some interesting results in that direction. 
For instance, a dilation theorem which derives from the work of Frazho 
\cite{Fra1}, Bunce \cite{Bun}, and Popescu \cite{Pop_diln} provides a 
generalization of Sz.-Nagy's classical minimal dilation of a contraction 
to an isometry \cite{SF}. Further, Popescu \cite{Pop_diln} gives an 
appropriate version of the classical Wold decomposition in this arena. In 
each of these cases, the role of unitary operators is played by $n$-tuples 
of isometries $S = \rows$ which satisfy (\ref{cuntz}). 
Moreover, the so-called `left creation operators' acting on Fock space Hilbert space, 
which also arise in theoretical physics and 
free probability theory, provide the appropriate generalization of 
unilateral shifts 
for this setting. They play a role in the aforementioned  dilation context, 
and the algebras they generate give noncommutative analogues of 
Toeplitz-type algebras. 
These algebras were discovered by Popescu \cite{Pop6,Pop_vN}, and there is now a growing
body of literature on these and other related algebras. 
For example, see work of Arias, Davidson, Katsoulis, 
Muhly, Pitts, 
Popescu, Power, Solel, and the author  
\cite{AP,Dsurvey,DKP,DP1,DP2,KP2,KP,Kribs_factor,
MS,Pop_fact,Pop_beur,Pop_func,Pop_disc}. 
Weighted versions of these shifts have been 
analyzed recently as well \cite{Aproject,AP2,nonsa,inductive}.

In this paper, we present and investigate versions $T = \rowt$ of 
bilateral weighted shift operators for the noncommutative 
multivariable setting.  
Specifically, our investigation will focus on the joint 
reducing subspace structure of the operators $T 
=\rowt$; or equivalently, reducibility of the weak operator 
topology 
closed algebras $\fsulam$ they generate. From one point of view we 
are considering certain ``weighted free semigroup algebras'', if we wish to 
keep in line with terminology of \cite{Dsurvey,DKP,DKS,DP1}. We derive a 
complete 
characterization of the reducing subspaces for $\fsulam$, and 
this 
gives a generalization of Nikolskii's Theorem for the 
single variable case \cite{Nikol}. 
The unweighted cases for these shifts, which we 
denote by $S = \rows$, may be 
regarded as generalizations of the canonical bilateral shift from the single variable 
theory. These $n$-tuples determine a subclass of  
the atomic representations of $\O_n$ considered by Davidson and Pitts \cite{DP1}.

\section{Noncommutative Bilateral Weighted Shifts}\label{S:shifts}

For succinctness we drop the `multivariable' reference when referring to these shifts. 
We first give a description of the possible Hilbert spaces. 
Throughout the rest of the paper $n\geq 2$ will be a fixed positive 
integer or $n= \infty$, but we will behave as though $n$ is finite. 

\begin{defn}\label{spacedefn}
Let $\fnplus$ be the unital free semigroup on $n$ noncommuting letters 
written as  
$\{1,2,\ldots, n\}$. Let $\Omega_+$ be the collection of infinite words in 
the generators of $\fnplus$. Given $\omega = i_1 i_2 \cdots $ in 
$\Omega_+$, define a sequence of words in $\fnplus$ by 
\[
\omega_m = \left\{ \begin{array}{ll}
i_1i_2\cdots i_m & \mbox{for $m\geq 1$} \\
\phi & \mbox{for $m=0$.}
\end{array}\right.
\]
The element $\phi$ is the unit in $\fnplus$, corresponding to the empty 
word. Let  $\fom$ denote the set of reduced words in the free group on $n$ 
generators of the form $u=v\omega_m^{-1}$ for $v\in\fnplus$ and some 
$m\geq 0$. Let $\K_\omega = \fock$ be the Hilbert space with orthonormal 
basis $\{\xi_u : u\in\fom\}$. 
\end{defn} 

For each Hilbert space $\K_\omega$, there is a class of bilateral weighted 
shifts which naturally act on it. 

\begin{defn}\label{shiftdefn}
An $n$-tuple of operators $R= (R_1,\ldots, R_n)$ acting on a Hilbert space 
$\H$ forms a {\it noncommutative bilateral weighted shift} if there is an 
$\omega\in\Omega_+$ and a unitary $U: \H \rightarrow \K_\omega$ for which 
there are scalars $\Lambda = \{\lambda_u\}_{u\in\fu}$ and operators $T = 
\rowt$, given by $T_i = U R_i U^*$, such that 
\begin{eqnarray}\label{shiftaction}
T_i \xi_u = \lambda_{iu} \xi_{iu} \qfor u\in\fu \qand 1\leq i \leq n.
\end{eqnarray}
We also define $\fsulam$ to be the unital weak operator topology closed 
(nonselfadjoint) algebra generated by $\{T_1, \ldots, T_n\}$,
\[
\fsulam = \wot \!\!-\!\!\Alg \{T_1,\ldots, T_n\}.
\] 
\end{defn}

\begin{notes}
{\bf ($\mathbf{i}$)} The shifts $T= \rowt$ we 
consider will be assumed to act on a space $\K_\omega$ as in (\ref{shiftaction}). 

\noindent{
{\bf ($\mathbf{ii}$)}} Given $\omega\in\Omega_+$, denote the operators 
for the unweighted case ($\lambda_u\equiv 1$) by $S= \rows$, and the 
associated algebra by $\fS_\omega$. Observe the operators $S = 
\rows$ determine a representation of $\O_n$. In fact, they form a subclass 
of the {\it atomic} representations of $\O_n$ classified by Davidson and 
Pitts \cite{DP1} in their ongoing investigation of {\it free semigroup 
algebras}, of which $\fS_\omega$ provides an example. 

\noindent{
{\bf ($\mathbf{iii}$)}} Every infinite word $\omega$ in $\Omega_+$ is 
either  
{\it aperiodic} or eventually periodic in the sense that there are 
$u,v_0\in\fnplus$ such that $\omega = uv_0v_0v_0\cdots$. But for our 
purposes, we lose no generality in focusing on the $\omega\in\Omega_+$ 
which are either aperiodic, or truly periodic ($u=\phi$). The reason is 
that given $\omega= u\omega^\prime$ in $\Omega_+$ with $u\in\fnplus$, 
$|u|=k$,   and 
$\omega^\prime \in\Omega_+$ periodic, one can define a unitary 
$U:\K_{\omega^\prime} 
\rightarrow \K_{\omega}$ by $U\xi_{v(\omega_m^\prime)^{-1}} = 
\xi_{v\omega_{m+k}^{-1}}$. This unitary 
intertwines the associated 
$n$-tuples $T^{(\omega)}$ and $T^{(\omega_0)}$ by $U^* T_i^{(\omega)}U = 
T_i^{(\omega^\prime)}$, and hence the algebras $U\fsulam U^* = 
\fS_{\omega^\prime,\Lambda}$. This type of equivalence is called a 
{\it 
shift-tail unitary equivalence}. 

\noindent{
{\bf ($\mathbf{iv}$)}} The shifts $T = \rowt$ acting on $\K_\omega$ 
trace out an infinite tree structure for $\K_\omega$, where basis vectors 
$\xi_u$ and weights $\lambda_u$ are identified with vertices in the tree. This point is illustrated
further in the examples of Section~\ref{S:examples}. The vertex set  
corresponding to  words $\{\omega_m^{-1}:m\geq 0\}$ is referred to as  the {\it 
main  branch} of $\fu$.  

\noindent{
{\bf ($\mathbf{v}$)}} We use the following notions of length for elements 
of $\fu$ and $\fnplus$: For $u_1, u_2\in\fu$, $|u_1|\leq |u_2|$ means that 
$u_1$ lies at least as close as $u_2$ to the closest common ancestor of 
$u_1,u_2$ on the main branch of $\fu$. Further, given $w\in\fnplus$, $|w|$ 
denotes the length of the word $w$. 

\noindent{
{\bf ($\mathbf{vi}$)}} Observe that for $n=1$, the set $\Omega_+$ consists 
of 
the single element $\omega = 111\cdots$. Hence $\K_\omega$ is a single 
two-way 
infinite stalk, and in Definition~\ref{shiftdefn} we recover standard 
bilateral weighted shift operators.  
\end{notes}

As we are interested in describing non-degenerative reducibility of $\fsulam$, we may clearly
assume the weights $\Lambda = \{\lambda_u \}_{u\in\fu}$ are nonzero. Furthermore, 
the result below shows we lose no generality in assuming weights are nonnegative. Thus the
following assumption on weights will be made throughout the paper:
\[
\mbox{{\bf Assumption:}} \,\,\,\,\, \lambda_u > 0 \qfor u\in\fu.
\]

\begin{prop}\label{unitequiv}
Let $T=\rowt$ be a bilateral weighted shift on $\K_\omega$ with weights 
$\Lambda 
= \{ \lambda_u \}_{u\in\bbF_\omega}$. Then there is a unitary $U$ in 
$\B(\K_\omega)$, 
which is diagonal with respect to the standard basis for $\K_\omega$, such 
that 
\[
(UT_1U^*, \ldots, UT_nU^*)
\]
is a bilateral weighted shift on $\K_\omega$ with weights $\Lambda^\prime 
= \{ 
|\lambda_u |\}_{u\in\bbF_\omega}$. Thus, the algebras $\fsulam$ 
and 
$\fS_{\omega,\Lambda^\prime}$ are unitarily equivalent. 
\end{prop}

\Prf
We define the unitary $U$ by $U\xi_v = \mu_v \xi_v$ upon inductively 
choosing scalars $\mu_v$ for $v\in\bbF_\omega$. Put $\mu_\phi = 1$. Then,  
assuming $\mu_{\omega^{-1}_{m-1}}$ has been chosen, for all  
$m\geq 1$ we  choose 
$\munegk\in\bbC$ of modulus one such that
\[
(\lambda_{\omega_{m-1}^{-1}} \mu_{\omega_{m-1}^{-1}}) \,\ol{\munegk} = 
c_{\omega_{m-1}^{-1}} 
\geq 0.
\]
Then we have $UT_{i_m} U^* \xi_{\omega_m^{-1}} = c_{\omega_{m-1}^{-1}} 
\xi_{\omega_{m-1}^{-1}}$ 
for $m\geq 1$, and this takes care of the main branch scalars. Now, if we 
are given a   word $v$
in $\bbF_\omega$ and $\mu_v$ has been chosen, then for $1\leq i \leq n$ 
with $iv$ off the main branch,  choose 
scalars $\mu_{iv}$ of modulus one such that 
$(\ol{\mu_v}\lambda_{iv})\mu_{iv}=c_{iv} \geq 0$. Thus, $UT_iU^* 
\xi_v = c_{iv} \xi_{iv}$ for all $v$ and $i$, and this yields the 
desired unitary. 
\bx
 
The following factorization result will be useful in the sequel. 

\begin{prop}\label{factor}
Let $T=\rowt$ be a bilateral weighted shift on $\K_\omega$. Then $T_i = 
S_i 
W_i$ 
for $1\leq i \leq n$, where $S=\rows$ is the unweighted shift on 
$\K_\omega$  and 
each $W_i$ is a positive operator, which is diagonal with respect to the 
standard basis for $\K_\omega$, given by 
\[
W_i \xi_u = \lambda_{iu} \xi_u 
\qfor 1\leq i \leq n \qand u\in\fu.
\]
Furthermore, a subspace $\H$ of $\K_\omega$ reduces the family  
$\{T_1, \ldots, T_n\}$, equivalently the algebra $\fsulam$, if and only if 
it reduces the family of operators $\{S_1, \ldots ,S_n, W_1, \ldots, W_n 
\}.$ 
\end{prop}

\Prf
The operators $T_i$ are easily seen to factor as $T_i=S_iW_i$, where $W_i$ 
is given by
$
W_i\xi_u = (T_i\xi_u, \xi_{iu}) \xi_u = \lambda_{iu} \xi_u,
$
for $u\in\fu$ and $1\leq i \leq n$. The last statement follows  from 
standard operator theory. Indeed, if $P$ is a projection which commutes 
with  $T_i=S_iW_i$, then it commutes with $ T_i^*T_i = W_i^2$, 
and hence with the 
positive operator $W_i= \sqrt{W_i^2}$. Thus, $PS_iW_i=S_iW_iP=S_iPW_i$ and 
$P$ commutes 
with $S_i$ on the range of $W_i$. But by our assumption on weights, $W_i$ 
is 
surjective, whence $S_iP=PS_i$. The other direction is obvious. 
\bx

\section{Main Results}\label{S:main}

In this section we investigate the reducing subspace 
structure of the algebras $\fsulam$ and, at the same time,  the joint 
reducing subspace structure of  weighted shifts $T = \rowt$. We begin by 
addressing the irreducible case. 

\begin{prop}\label{aperiodic}
Let $\omega$ be an aperiodic word in $\Omega_+$. Then the algebra 
$\fsulam$ is 
irreducible. 
\end{prop}

\Prf
First note that the projections $\{\omprojns : m\geq 0\}$ belong to 
the von Neumann algebra $\fA_T$ generated by $\{T_1, \ldots, T_n\}$. This 
follows from Proposition~\ref{factor}, because the von Neumann algebra 
$\fA_S$ generated by the unweighted shifts $\{S_1,\ldots, S_n\}$ is 
spanned by its projections, and hence the double commutant identity 
implies 
$\fA_S = \fA_S^{\prime\prime} \subseteq \fA_T^{\prime\prime} = \fA_T$. 

But when $\omega$ is an aperiodic word, the projections $\omprojns$  
converge in the strong operator topology to the rank one projection 
$\xi_\phi \xi_\phi^*$ onto the subspace $\spn\{\xi_\phi\}$. Indeed, the vector 
$\xi_\phi$ clearly belongs to the range of the projection 
$P=\sotlim_{m\rightarrow\infty} \omprojns$, which exists since the 
$\omprojns$ form a decreasing sequence of projections. On the other hand, 
given a word $u=v\omega_k^{-1}\in\fu$ and positive integer $m$, there is a 
unique word $u_m\in\fnplus$ of length $m$ such that $\xi_u$ is in the 
range of 
$u_m(S)$. This determines an infinite word in $\Omega_+$, which 
is $v(\omega_k^{-1}\omega)=u\omega$, and $u_m$ makes up the first $m$ 
terms of this word. 
As $\omega$ is aperiodic, $u\omega$ and $\omega$ are different except 
when $u=\phi$; for if $u\omega = \omega$, then $|v|\neq k$ implies 
$\omega$ is periodic, whereas $|v|=k$ implies $v=\omega_k$ and $u=\phi$.   Hence for 
some sufficiently large $m$, the words $u_m$ 
and $\omega_m$ are distinct when $u\neq \phi$. Whence, $\omprojns \xi_u = 0$, and it follows
that these
projections converge to $\xi_\phi \xi_\phi^*$.   

Thus, it follows that a projection $Q\in\fA_T^\prime$ commuting with  
 the operators $\{T_1,\ldots, T_n\}$,  
also commutes with $P= \xi_\phi \xi_\phi^*\in \fA_T$. But the vector $\xi_\phi$ 
is clearly cyclic for $\fA_T$ since all weights are assumed to be nonzero. Hence 
$Q=0$ or $Q=I$, as required.  
\bx

\begin{rem}
We mention that the algebras $\fsulam$ coming from aperiodic words 
$\omega$ 
considered in Theorem~\ref{aperiodic} truly are exclusive to the 
noncommutative setting. Indeed, the only infinite word $\omega=111\cdots$ 
with 
letters in a one letter alphabet is, of course, trivially periodic. 
\end{rem}

We turn now  to the periodic case. 
Given a word $\omega = v_0v_0v_0\cdots $ in $\Omega_+$  with 
$v_0\in\fnplus$ {\it primitive} (not a power of another element in 
$\fnplus$) and a 
positive integer $k\geq 1$, there is a canonical decomposition of $\fu$ 
into mutually disjoint subsets, and hence of the Hilbert space $\K_\omega 
= \ell^2(\fu)$  into mutually orthogonal subspaces. First define for $l \leq 0$
\[
\bbF_{\omega,k}^{(l)} = \big\{ v\,\omega_m^{-1} :v\in\fnplus\setminus \fnplus i_m,\, -lk|v_0|
\leq m < (-l+1)k|v_0| \big\},
\]
where we take $i_0 = i_{|v_0|}$. For $l \geq 1$ define 
\[
\bbF_{\omega,k}^{(l)} = \big\{ v\,\omega^\prime_{m+1} :v\in\fnplus\setminus \fnplus
i^\prime_m,\,
(l-1)k|v_0| < m \leq  l k|v_0| \big\},
\]
where $\omega^\prime = \cdots v_0v_0$, and $\omega_m^\prime\in\fnplus$ is the initial
segment (from the right) of $m$ letters in this infinite word, with $i_m^\prime$ equal to the
$m$th letter in the
sequence. We refer 
to $\bbF_{\omega,k}^{(0)}:=\fomk$ as the {\it principal component} of this 
partition of $\fu$. Observe that the unit $\phi$ belongs to $\fomk$.  Given $u\in\fomk$, 
further define $\ju$ to be the subset of $\fu$ 
consisting of $u$ and all its natural translates in the sets $\{ \fomkl : l\in\bbZ\} $. 
Specifically, define 
\[
\J_{\phi,k} = \big\{ (v_0^k)^l : l\in\bbZ \big\} \qand \ju = u\J_{\phi,k} \qfor u\in\fomk. 
\] 
Thus, in summary, given a periodic word $\omega\in\Omega_+$ and a positive integer $k\geq
1$, there is a canonical partition of $\fu$ into a disjoint union 
of  subsets which generates a spatial decomposition of $\K_\omega$; 
\[ 
\fu = \bigcup_{u\in\fomk} \ju \qand \K_\omega = \sum_{u\in\fomk}\! 
\!\oplus 
\,\spn\{ \xi_v:v\in\ju\}.
\] 

\begin{defn}\label{periodicdefn}
Let $\omega= v_0v_0\cdots$ be a periodic word in $\Omega_+$ and let 
$k\geq 1$ be a positive integer. Let $\fu = \cup_{u\in\fomk} \ju$ be the 
associated canonical decomposition of $\fu$. We shall say that a bilateral 
weighted shift $T= \rowt$ is {\it period} $k$ (or equivalently, 
the 
weights $\Lambda= \{\lambda_u\}_{u\in\fu}$ are period $k$) if 
\[
T_i\, \xi_v = \lambda_{iv}\,\xi_{iv} = \lambda_u\, \xi_{iv} \qfor v\in\fu 
\qand 
1 \leq i \leq n,
\]
where $u$ is the unique element of $\fomk$ with $iv\in\ju$. 
\end{defn}

\begin{note}
The weights $\{\lambda_u: u\in\fomk \}$ may be regarded as the 
{\it remainders} of a $k$-periodic shift. There is a satisfying visual 
interpretation of this periodicity in terms of the infinite tree structure 
of $\K_\omega$. This is expanded on in the examples of the next section. We 
discovered this notion of periodicity while proving the following 
theorem.
\end{note}

\begin{thm}\label{periodic}
The following assertions are equivalent for a periodic word $\omega$ in $\Omega_+$: 
\begin{itemize}
\item[$(i)$] $\fsulam$ has non-trivial reducing subspaces. 
\item[$(ii)$] The weights $\Lambda = \{\lambda_u\}_{u\in\bbF_\omega}$ are 
periodic. 
\item[$(iii)$] The weighted shift $T= \rowt$ is periodic.
\end{itemize}
\end{thm}

\Prf
The last two conditions are equivalent by definition. The 
implication $(ii)\Rightarrow (i)$ is established in the proof of 
Theorem~\ref{mainthm} below. 
We prove $(i)\Rightarrow (ii)$. Suppose $\H$ is a non-trivial subspace 
of $\K_\omega$ which reduces $\fsulam$. Let $P_\H \notin \{0,I\}$ be the 
projection onto $\H$. Then  $P_\H T_i = T_i P_\H$ for 
$1 \leq i \leq n$, and hence $P_\H$ commutes 
with the family $\{S_i, W_j: 1 \leq i,j \leq n\}$ by 
Proposition~\ref{factor}. 

Consider the set $\P$ of all partitions $\Pi$ of the set $\fu$ 
into 
mutually disjoint subsets such that
\[
Q_\J P_\H = P_\H Q_\J \qfor \J \in \Pi,
\]
where $Q_\J$ is the projection onto $\spn\{\xi_u:u\in\J\}$. 
Notice that  $\P$ is closed under a natural join operation; if 
$\Pi_1, \Pi_2\in\P$, then $\Pi_1\bigvee\Pi_2 = 
\{\J : \J=\J_1\cap\J_2,\J_i\in\Pi_i,i=1,2\}$ belongs to $\P$. This set is 
clearly  a partition of $\fu$ into disjoint subsets, and the corresponding projections $Q_\J =
Q_{\J_1\cap\J_2} = Q_{\J_1}Q_{\J_2}$ commute with $P_\H$. 

The following notation will be useful. 
If we are given $\J\subseteq \fu$ such that $Q_\J P_\H = P_\H 
Q_\J$, then for all words $w\in\fnplus$ define subsets of $\fu$ by  
\[
w(S) \J \equiv \{ wu: u\in\J \} = w\J \qand w(S)^*\J \equiv w^{-1}\J\cap \fu,
\]
where $w(S)$ is the isometry $w(S) = S_{i_1}\cdots S_{i_k}$ when $w=i_1\cdots i_k$. 
Observe that the projections $Q_{w(S)\J}$ and  $Q_{w(S)^*\J}$ are, respectively, the
projections onto the ranges of the operators $w(S)Q_\J$ and $w(S)^*Q_\J$. 
But $w(S)Q_\J$ is a partial isometry, so that 
\[
Q_{w(S)\J } = ( w(S) Q_\J)  ( w(S) Q_\J)^* = w(S) Q_\J w(S)^*. 
\]
Moreover, the projection $Q_\J$ leaves the range of $w(S)$ invariant since the ranges of both
these operators are  spanned by standard basis vectors. Thus, $w(S)^* Q_\J$ is also a partial
isometry with final projection 
\begin{eqnarray*}
Q_{w(S)^*\J} &=& ( w(S)^*\, Q_\J\, w(S) ) ( w(S)^* \,Q_\J\, w(S) ) \\
&=&  w(S)^* \,Q_\J\, (w(S) w(S)^*)\, Q_\J\, w(S) =  w(S)^*\, Q_\J\, w(S).  
\end{eqnarray*} 
Hence these projections commute with $P_\H$ as well. Further  
define subsets $v(S)w(S)^*\J$, with $v,w\in\fnplus$, in a similar manner, yielding  projections 
$Q_{v(S)w(S)^*\J}= v(S) w(S)^*\, Q_\J\, w(S) v(S)^*$ that commute with $P_\H$ as well.  

Next define a distinguished partition $\Pi_0$ in $\P$, determined by the equivalence 
relation 
\[
u\sim v \qfor u,v\in\fu \qiff \lambda_u = \lambda_v.
\] 
To see that $\Pi_0$ belongs to $\P$, for $r> 0$ let $Q_r$ be the 
projection onto $\spn\{ \xi_u: \lambda_u = r\}$. If nonzero, this is a 
typical projection onto a subspace of $\K_\omega$ determined by a coset of 
$\Pi_0$. 
For $1\leq i \leq n$, let $\pri$ be the projection onto the eigenspace
\[
\pri\K_\omega = \ker\big(  W_i-rI \big) = \spn\{ \xi_u : \lambda_{iu}=r\}. 
\]
As spectral projections for the $W_i$, the $P_{r,i}$ commute with 
$P_\H$ and so do the operators $S_i\pri$. Let $\qri = S_i\pri S_i^*$ be the projection 
onto the range of the partial isometry $S_i\pri$, which is $\spn\{ \xi_v: \lambda_v = r, v=iu,
u\in\fu \}$. Then $Q_r
= \sum_i \oplus \,\qri$ is a 
projection which evidently commutes with $P_\H$. 
Thus $\Pi_0$ belongs 
to $\P$ as claimed.

We define a partial ordering on the partitions $\P$ by: $\Pi_1\geq \Pi_2$ 
if every $\J_2\in\Pi_2$ may be obtained by $\J_2 = 
\bigcup_{\J_1\subseteq\J_2}\J_1$ with $\J_1\in\Pi_1$. Let $\P_0= \{ 
\Pi\in\P : \Pi \geq \Pi_0 \}$. We claim that every chain $\{ \Pi_\alpha : 
\alpha \in \bbA\}$ in $\P_0$ has a maximal element $\Pi^\prime$ in $\P_0$. 
Indeed, define the partition $\Pi^\prime$ by the equivalence relation 
$u\sim v$, for $u,v\in\fu$, if for any $\alpha\in\bbA$, there exists 
$\J_\alpha\in\Pi_\alpha$ such that $u,v\in\J_\alpha$. 
This is an equivalence relation since each 
$\Pi_\alpha$ is a partition of $\fu$ into disjoint subsets. As the $\Pi_\alpha$ form a chain, it is
clear that  $\Pi^\prime\geq \Pi_\alpha \geq \Pi_0$ for all $\alpha\in\bbA$. Let 
$\J\in\Pi^\prime$. Observe from the definition of this relation, that for all $\alpha\in\bbA$, there
is a unique $\J_\alpha\in\Pi_\alpha$ with $\J \subseteq \J_\alpha$. On the other hand, if $u\in
\cap_{\alpha\in\bbA} \J_\alpha \supseteq \J$, then $u\sim v$ for all $v\in\J$. Hence $u\in\J$, and
we have $\J = \cap_{\alpha\in\bbA} \J_\alpha$.  Further, we have 
$\J_\alpha \supseteq \J_\beta \supseteq \J$ for $\Pi_\alpha \leq 
\Pi_\beta$ by the definition of the ordering. Thus, it follows that $P_\H$ 
commutes with the projection $Q_\J = Q_{\cap_\alpha\J_\alpha} = 
\bigwedge_\alpha Q_{\J_\alpha}$. Hence $\Pi^\prime$ belongs to $\P_0$ and 
majorizes the chain $\{ \Pi_\alpha : \alpha\in\bbA\}$ as claimed. 

Now apply Zorn's Lemma to obtain a maximal element $\ol{\Pi}$ of  
$\P_0$. We set aside a pair of  technical results on $\ol{\Pi}$ which 
we require:

\begin{lem}\label{minlprojns}
Let $\J\in\ol{\Pi}$. For $m\geq 0$, there is a unique word $w\in\fnplus$, 
$|w|=m$, such that $w(S)^*\J\neq \emptyset$. It follows that  the 
projections  
$
\big\{ Q_{v(S)w(S)^*\J} : v\in\fnplus \big\}
$
are minimal projections with ranges spanned by standard basis vectors in the commutant
$\{P_\H\}^\prime$. In particular,
if $\J_1\in\ol{\Pi}$ is such that $v(S)w(S)^*\J \cap \J_1 \neq \emptyset$, 
then $v(S)w(S)^*\J = \J_1$.  
\end{lem}

\Prf
Given $m\geq 0$, there is some word $w\in\fnplus$, $|w|=m$, for which 
$w(S)^*\J\neq \emptyset$ since $\J\neq \emptyset$. As 
$w(S)w(S)^*\J\neq\emptyset$, we in fact have $w(S)w(S)^*\J = \J$. Indeed, 
we clearly have $w(S)w(S)^*\J \subseteq \J$, and if this were a strict 
inclusion we could refine 
$\ol{\Pi}$ in $\P_0$ by replacing $\J$ by $w(S)w(S)^*\J$ and $\J\setminus 
w(S)w(S)^*\J$. It follows that there can be no other word 
$w^\prime\in\fnplus$, 
$|w^\prime|=m$, with $w^\prime(S)^*\J\neq \emptyset$.

Suppose $Q\leq Q_{v(S)w(S)^*\J} = v(S)w(S)^*\, Q_\J\,  w(S)v(S)^*$ is a nonzero projection in 
$\{P_\H\}^\prime$ with range spanned by standard basis vectors. Then it follows that
$w(S)v(S)^*\,Q \,v(S) w(S)^* $ is a nonzero projection in $\{P_\H\}^\prime$, whose range is
spanned by standard basis vectors, with 
\begin{eqnarray*}
 w(S)v(S)^*\,\,Q \,\,v(S) w(S)^* &\leq & w(S) w(S)^*\,\, Q_\J \,\,w(S) w(S)^* \\               
&=& Q_{w(S) w(S)^*\J} = Q_\J.  
\end{eqnarray*}
Hence $w(S)v(S)^*\, Q\, v(S) w(S)^* =Q_\J$ 
by the minimality of $Q_\J$ as a projection spanned by standard basis vectors in
$\{P_\H\}^\prime$ (which 
follows from the maximality of the partition $\ol{\Pi}$). Thus it follows  that 
\begin{eqnarray*}
Q &=& v(S)v(S)^*\,\,Q\,\, v(S) v(S)^* \\
 &=&  v(S)w(S)^* \big( w(S)v(S)^* \,\, Q \,\,v(S) w(S)^* \big) w(S) v(S)^* \\
& =&  Q_{v(S)w(S)^*\J}.
\end{eqnarray*}

Finally, let $\J_1\in\ol{\Pi}$ be such that $v(S)w(S)^*\J \cap \J_1\neq \emptyset$. 
Then $\J_1\subseteq v(S)w(S)^*\J$; otherwise we could refine $\ol{\Pi}$ 
by replacing $\J_1$ by $\J_1\cap v(S)w(S)^*\J$ and $\J_1\setminus 
v(S)w(S)^*\J$. But by the previous minimality argument we in fact have 
$\J_1 = v(S) w(S)^*\J$, as required. 
\bx

\begin{lem}\label{mainbranchelts}
Every $\J\in\ol{\Pi}$ contains more than one element. (In fact, we will 
see that every $\J\in\ol{\Pi}$ contains infinitely many elements.) 
Moreover, there is a 
$\J_1\in\ol{\Pi}$ with two elements lying on the same path in $\fu$, and it follows 
that there is a word $w\in\fnplus$ such that $w(S) \J_1= \J_1$. In 
particular, $\J_1$ contains elements lying on the main branch of $\fu$. 
\end{lem}

\Prf
Let $\J\in\ol{\Pi}$. Then $\J$ is not a singleton. Suppose to 
the contrary that $\J = \{u\}$. For all $v\in\fu$, there are words $v_+$, 
$v_-$ in $\fnplus$ such that
$
v_+(S) v_-(S)^* \xi_u = \xi_v.
$
In particular, the nonempty subsets amongst $\{ v_+(S)v_-(S)^* \J^\prime :\J^\prime \in
\ol{\Pi}\}$ are part of a partition in $\P$ which includes the 
subset $\{v\}$. Since $v\in\fu$ was arbitrary, and $\P$ is closed 
under the meet operation, it follows that the partition 
$\big\{\{v\}:v\in\fu\big\}$ belongs to $\P$. Hence each vector $\xi_v$ 
either belongs to $\H$ or is orthogonal to it. Thus, as $\H\neq\{0\}$, 
it includes some standard basis vector $\xi_{v_0}$. However, every standard basis vector  is
cyclic for the von Neumann algebra $\fA_S$, and $P_\H$ belongs to the commutant
$\fA_S^\prime$. This implies, incorrectly, that 
$\H = \K_\omega$ since $\H \supseteq \ol{\fA_S \xi_{v_0}} = \K_\omega$. Thus $\J$ contains
at 
least two elements. 

Now suppose $u_1\neq u_2$ belong to $\J$ with $|u_2|\geq |u_1|$. Let 
$u\in\fu$ be the nearest 
common ancestor of $u_1,u_2$. Then $u$ lies closer to $u_1$ in the tree 
structure for
$\K_\omega$. Suppose 
$u\in\J_1\in\ol{\Pi}$. Choose $v_1\in\fnplus$ such that $\xi_{u_1} = 
v_1(S)\xi_u$. Then
$v_1(S)\J_1 \cap \J \neq \emptyset$, and hence we  have $\J = 
v_1(S)\J_1$ by Lemma~\ref{minlprojns}. Thus $\xi_{u_2} = v_1(S)\xi_x$ for 
some $x\in\J_1$. Since $|u_2|\geq |u_1|$, we have $|x|\geq |u|$, and the  
existence of a $w\in\fnplus$ such that $w(S)\xi_u = \xi_x$ follows from 
$u$ and $x$ both being ancestors of $u_2$. As 
$u=x$ would 
incorrectly imply $\xi_{u_1}= \xi_{u_2}$, we have shown that $\J_1$ 
contains elements $u\neq x$ for which $w(S) \xi_u = \xi_x$ for some 
$w\in\fnplus$. But another application of Lemma~\ref{minlprojns} yields 
$w(S)\J_1 = \J_1$, since these sets have non-trivial intersection. Thus 
$\big(w(S)^*\big)^m\J_1 = \J_1$ for $m\geq 0$. It follows that $\J_1$ 
contains a pair of elements (infinitely many in fact) which lie on the 
main branch. 
\bx

Let $\J_1\in\ol{\Pi}$ be obtained as in Lemma~\ref{mainbranchelts}, and 
let $u\in\fnplus$ be a word of minimal length such that $u(S)\J_1=\J_1$. 
We may assume with no loss of generality that the unit $\phi\in\J_1$. Indeed,  
there is a $v\in\fnplus$ with $\phi\in 
v(S)\J_1\in\ol{\Pi}$, and evidently $v(S)\J_1$ will satisfy the conditions 
on $\J_1$ 
in Lemma~\ref{mainbranchelts}.  
Thus the vectors  $(u(S)^*)^m\xi_\phi $, $m\geq 0$, are standard basis vectors on the main 
branch. It follows that if $\omega = v_0v_0v_0\cdots$ with $v_0\in 
\fnplus$ a primitive word, then there exists positive integers $m_1$, 
$m_2$ for which 
$v_0^{m_1} = u^{m_2}$. Thus we can apply the following combinatorial 
lemma: 

\begin{lem}\label{combinatoric}
Suppose $u, v\in \fnplus$ and $l, m\geq 1$ are positive integers 
such that $u^l = v^m$. If $v$ is a primitive word, then $u =v^k$ for some 
$k\geq 1$. 
\end{lem}

\Prf
Suppose first that $|v|\leq |u|$. We can clearly assume  $l,m>\!>0$. Then 
$u^l = v^m$ implies the 
existence of $u_1,u_2\in\fnplus$ such that $u=vu_1$ and $u=u_2v$ with 
$|u_1|=|u_2|$.  
If $u_1,u_2 \neq\phi$ (or $u\neq v$), then $u_2=vu_3$, $u_1= 
u_4v$, but 
$vu_3v=u=vu_4v$, whence $u_3=u_4$. If $u_3\neq \phi$ (or $u\neq v^2$), 
then $u_3=vu_5=u_6v$ with $|u_5 | = |u_6|$. If $u_5,u_6 \neq \phi $ (or 
$u\neq v^3$), then $u_5=u_7v$, $u_6= vu_8$ and $u_7=u_8$. We can iterate  
this process just finitely many times, hence we eventually get $u=v^k$ for 
some $k \geq 1$, as required. 
On the other hand, by the previous argument the $|u|\leq |v|$ case forces 
$u=v$ by primitivity of $v$.
\bx

Therefore, we have $u=v_0^k$ for some $k\geq 1$, and it follows that
\[
\J_1 = \big\{ (v_0^k)^l: l\in\bbZ \big\}\equiv \J_{\phi,k}. 
\]
Indeed, this set is contained in $\J_1$ as $\phi\in \J_1  = u(S)^l \J_1 = 
(u(S)^*)^l \J_1$, for $l \geq 0$. Also, the minimality of the length of $u$ such that $u(S)\J_1 =
\J_1$ ensures that the elements $\{ (v_0^k)^l : l \leq 0\}$ form the intersection of $\J_1$ and the
main branch. Thus 
if $v\in\J_1$, then some power $(u(S)^*)^l \xi_v = \xi_{u^{-l}v}$ will 
reside on the main branch, forcing $u^{-l}v$, and hence $v= u^l 
(u^{-l}v)$ to belong to this set. 

We have just shown that the set $\J_1$ coincides with the set $\J_{\phi,k}$ 
from the discussion prior to Definition~\ref{periodicdefn}. Thus, it follows from our analysis of
the partition $\ol{\Pi}$ that 
\[
\ol{\Pi} = \big\{  u \J_{\phi,k} : u\in \fomk \big\} 
= \big\{ \ju : u\in\fomk \big\}. 
\]
Therefore, the periodicity of $\fsulam$ and $T=\rowt$ is now an immediate consequence of   
$\ol{\Pi}\geq \Pi_0$. This completes the proof of 
Theorem~\ref{periodic}.
\bx

\begin{rem}
This result generalizes Nikolskii's Theorem for the 
single variable case (Theorem~4 from \cite{Nikol}), which asserts that 
bilateral weighted shift operators on Hilbert space have non-trivial 
reducing subspaces if and only if the associated weight sequence is 
periodic. In 
fact, the reader of \cite{Nikol} will notice that we have been inspired by 
the  proof of this theorem in establishing the implication 
$(i)\Rightarrow(ii)$ of Theorem~\ref{periodic}. However, we wish to emphasize 
that, evidently, there are a number of technical details here which do not  
arise in the single variable case.  
\end{rem}

We now combine the previous results and show how, in the doubly-periodic case,  non-trivial
reducing subspaces are determined by reducing subspaces of the free semigroup algebra
$\fS_\omega$. 
We require another  decomposition of $\fu$ when $\omega$ is periodic (see 
the discussion prior to Definition~\ref{periodicdefn}). 
For $l\in \bbZ$, partition $\fomkl$ into $k$ disjoint subsets $\fomkl= \bbG_0^{(l)} \cup \ldots
\cup   \bbG_{k-1}^{(l)}$. To avoid cumbersome notation, we only explicitly define the
decomposition
of the principal component $\fomk \equiv \bbF_{\omega,k}^{(0)} = \bbG_0^{(0)} \cup \ldots
\cup   \bbG_{k-1}^{(0)}$, the other cases are similar. Let 
\[
\bbG_r^{(0)} = \big\{ v\omega_m^{-1} : v\in\fnplus, r |v_0| \leq m < (r+1)  |v_0| \big\} \qfor
0\leq r < k. 
\]
Let $\K_{\omega,k}$ be the subspace of $\K_\omega$ defined by  
\[
\K_{\omega,k} := \sum_{m\in\bbZ} \oplus \spn \{ \xi_u : 
u\in\bbG_0^{(mk)} \}.
\]
Given $m\in\bbZ$, there is a natural bijection between the sets 
$\bbG_0^{(mk)}$ and $\bbG_0^{(m)}$ which determines a unitary operator
$\uomk : \K_{\omega,k} \rightarrow \K_\omega$ mapping $\K_{\omega,k}$ onto 
$\K_\omega = 
\sum_{m\in\bbZ}^\oplus 
\spn\{  \xi_u :u\in\bbG_0^{(m)} \}$.

\begin{thm}\label{mainthm}
The  following assertions are equivalent for an infinite word $\omega$ in $\Omega_+$:
\begin{itemize}
\item[$(i)$] The algebra $\fsulam$ has non-trivial reducing subspaces. 
\item[$(ii)$] The operators $T =\rowt$ have non-trivial joint reducing 
subspaces. 
\item[$(iii)$] The pair $(\omega,\Lambda)$ is doubly-periodic; that is, the word $\omega$ is
periodic and the weights $\Lambda=\{\lambda_u\}_{u\in\fu}$ are periodic. 
\end{itemize}
Furthermore, in the case that the word $\omega = v_0v_0\cdots$ is periodic 
and the weights $\Lambda=\{\lambda_u\}_{u\in\fu}$ are of period $k$, the 
reducing subspaces $\H\subseteq\K_\omega$ for $\fsulam$, equivalently the 
joint reducing subspaces for $T=\rowt$, are generated by the reducing subspaces for the free
semigroup algebra $\fS_\omega$ in the following manner:
\begin{eqnarray}\label{reduceform}
\H = \sum_{u\in\fomk} \oplus\,\, Q_{\ju}\H =\sum_{r=0}^{k-1} \oplus 
\,\, P_r \H, 
\end{eqnarray}
where the projections 
\begin{eqnarray}\label{projndecomp}
P_r = \sum_{u\in\bbG_r^{(0)}} \oplus \,\,Q_{\ju} =  \sum_{u\in\bbG_r^{(0)}} \oplus
\,\,u(S)Q_{\J_{\phi,k}} u(S)^*,
\end{eqnarray}
for $0 \leq r < k$, commute with $P_\H$. The subspace $P_0\H$ is contained in 
$\K_{\omega,k}$, and the image $\uomk P_0 \H$ is a reducing subspace for 
the free 
semigroup algebra $\fS_\omega$. Conversely, every reducing 
subspace $\H^\prime$ for the free semigroup algebra $\fS_\omega$ generates 
a reducing subspace $\H$ for $\fsulam$ in this way through the equation 
$P_0\H \equiv \uomk^*\H^\prime$. 
\end{thm}

\Prf
The only thing left to establish is the form of the reducing subspaces for 
$\fsulam$ when both $\omega$ and $\Lambda$ are periodic. 
The form of $\H$ in (\ref{reduceform}) follows from the proof of 
Theorem~\ref{periodic}. In that proof, we saw how the projections 
$Q_{\ju}$ satisfy $Q_{\ju} = Q_{u\J_{\phi,k}} = u(S) Q_{\J_{\phi,k}}u(S)^*$. 
The fact that 
$\H^\prime = \uomk P_0\H$ reduces $\fS_\omega$, equivalently $S_1, \ldots 
,S_n$, follows from; the $k$-periodicity of the weights $\Lambda$, equations
(\ref{reduceform}) and  (\ref{projndecomp}), and the projections 
$Q_{\J_{u,k}}$ commuting with $P_\H$.   We 
leave this technical detail to the interested reader. 
On the other hand, every reducing subspace $\H^\prime$ for $\fS_\omega$ 
can be seen to generate an $\fsulam$-reducing subspace by first defining $P_0\H 
\equiv \uomk^*\H^\prime$, and then obtaining the rest of the subspaces $P_r\H$ by translations.  
\bx

\begin{rem}\label{unweightedsubspaces}
For the sake of brevity, we have deliberately avoided some technical 
details in giving the characterization of reducing subspaces for $\fsulam$ 
in the doubly-periodic case of Theorem~\ref{mainthm}. Suffice it to say, 
the reducing subspaces for 
$\fsulam$ are precisely those subspaces which are 
generated in a natural way by reducing subspaces for the unweighted 
shifts, 
in a manner which is 
analogous to the single variable case $n=1$ \cite{Nikol}. We 
mention that  reducing subspaces for the unweighted shifts $S= \rows$, 
equivalently for
$\fS_\omega$, were discussed by Davidson and Pitts \cite{DP1} in the context  of
representation theory for free semigroup algebras, hence we shall not pursue 
this here.  
\end{rem}

\section{Periodic Examples}\label{S:examples}

In this section, we explore further the notion of periodicity discovered in
Theorem~\ref{periodic}.  In particular, by considering the simplest possible examples, we show
there are natural infinite `shift matrices' associated with periodic $n$-tuples $T = \rowt$. This
gives information on the various operator algebras generated by these shifts, and leads us to a
number of open problems posed at the end of this section. 

\begin{eg}\label{example1}
Let $n=2$ and let $\omega = 222\cdots$ belong to $\Omega_+$. We shall consider weighted
shifts $T =(T_1,T_2)$ acting on $\K_\omega$ of period $k=1$ and $k=2$. Unlike the single
variable case, there is  diversity
even amongst the class of 1-periodic shifts. 

\vspace{0.05in}

{\noindent}${\mathbf (i)}$ Let $T =(T_1,T_2)$ be a weighted shift on 
$\K_\omega$ of period $k=1$. The
1-periodicity of $T$ implies the actions of $T_1,T_2$ are completely determined by their
behaviour on the principal component $\bbF_{\omega,1}$, which has the following natural
lexicographic ordering, 
\begin{eqnarray*}
\bbF_{\omega,1} &=& \big\{ \{\phi\}, \{v\in\fnplus 1\} \big\} \\
&=& \big\{ \{\phi\}, \{1\}, \{1^2, 21\}, \{1^3, 121, 21^2, 2^21\}, \ldots \big\}.
\end{eqnarray*}
The Hilbert space $\K_\omega$ has an infinite tree structure traced out by the operators
$T_1$, $T_2$. Vertices are identified with standard basis vectors $\{ \xi_u : u\in\fu\}$ and
weights $\Lambda = \{\lambda_u\}_{u\in\fu}$. The tree structure is as follows in this
case:

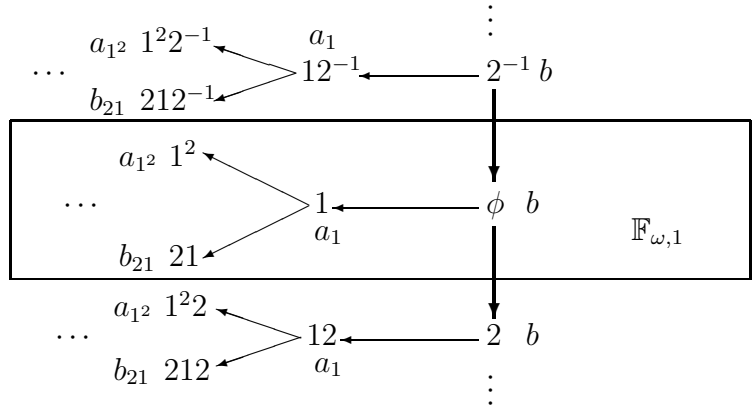
\begin{picture}(100,200)(-90,0)


\put(100,50){$2$}
\put(115,50){$b$}
\put(100,100){$\phi$}
\put(115,100){$b$}
\put(100,150){$2^{-1}$}
\put(120,150){$b$}
\multiput(100,165)(0,5){3}{$\cdot$}
\multiput(100,35)(0,-5){3}{$\cdot$}

\put(155,90){$\bbF_{\omega,1}$}

{\thicklines      
\put(103,95){\vector(0,-1){35}}
\put(103,147){\vector(0,-1){35}}
}

\put(-80,135){\line(1,0){280}}
\put(-80,75){\line(1,0){280}}
\put(-80,75){\line(0,1){60}}
\put(200,75){\line(0,1){60}}

\put(35,100){$1$}
\put(35,90){$a_1$}
\put(97,102){\vector(-1,0){55}}

\put(32,50){$12$}
\put(35,40){$a_1$}
\put(97,52){\vector(-1,0){52}}

\put(30,150){$12^{-1}$}
\put(33,164){$a_1$}
\put(97,152){\vector(-1,0){45}}

\put(33,103){\vector(-2,1){40}}
\put(33,103){\vector(-2,-1){40}}
\put(-20,118){$1^2$}
\put(-20,80){$21$} 
\put(-39,80){$b_{21}$}
\put(-39,118){$a_{1^2}$} 

\multiput(-60,100)(5,0){3}{$\cdot$}

\put(30,53){\vector(-3,1){32}}
\put(30,53){\vector(-3,-1){32}}
\put(-22,61){$1^22$}
\put(-22,37){$212$} 
\put(-41,37){$b_{21}$}
\put(-41,61){$a_{1^2}$} 
\multiput(-63,50)(5,0){3}{$\cdot$}

\put(28,153){\vector(-3,1){31}}
\put(28,153){\vector(-3,-1){31}}
\put(-31,161){$1^2 2^{-1}$}
\put(-31,139){$212^{-1}$} 
\put(-50,139){$b_{21}$}
\put(-50,161){$a_{1^2}$} 
\multiput(-72,150)(5,0){3}{$\cdot$}

\end{picture}

\noindent{As} an illustration of how $T = (T_1,T_2)$ acts on basis vectors, we note  
\[
\left\{ \begin{array}{rcl}
T_1\xi_\phi &=& a_1 \xi_1 \\
 T_2 \xi_\phi &=& b \,\xi_2
\end{array}\right. 
\,\,\,\,\, {\rm and} \,\,\,\,\,
\left\{ \begin{array}{rcl}
T_1\xi_{1} &=& a_{1^2} \xi_{1^2} \\
 T_2 \xi_{1} &=& b_{21} \xi_{21} 
\end{array}\right. 
\]

With our previous notation, the set $\J_{\phi,1}$ here consists of all the 
main branch elements together with the positive powers,  
$\J_{\phi,1} = \{ 2^l : l\in\bbZ\}$, and $\J_{u,1} = u \J_{\phi,1}$ for $u\in\bbF_{\omega, 1}$. 
Consider the subspaces 
\[
\H_u = \spn \{ \xi_v : v\in\J_{u,1} \} = u(S) \H_\phi \qfor u\in\bbF_{\omega,1}. 
\]
Observe that $\K_\omega = \sum_{u\in\bbF_{\omega,1}}^\oplus \H_u$ since 
$\fu$ is partitioned by $\{ \J_{u,1} : u\in\bbF_{\omega,1} \}$, and the 
adjoint operator 
$u(S)^*$ acts as a unitary from $\H_u = u(S) \H_\phi$ onto $\H_\phi$. Thus, putting these
equivalences together yields $\K_\omega \simeq \sum_{u\in\bbF_{\omega,1}}^\oplus \H_\phi$,
and if we order the index set $\bbF_{\omega,1}$ as above, then we obtain the following block
matrix decompositions for $T_1$ and $T_2$ (up to unitary equivalence):
\[
T_1 \simeq \left[
\begin{matrix}
\begin{array}{cc|cc|}
0 & 0 & 0 & 0 \\
a_1 I  & 0 & 0 & 0 \\ 
\hline
0 & a_{1^2}I & 0 & 0 \\
0&0&0&0 \\ 
\hline
0&0& a_{1^3}I & 0 \\
0&0&0& a_{121}I  \\
0&0&0&0 \\
0&0&0&0 \\
\hline
\end{array} 
  &  \cdot & \cdot & \cdot & & \\

\cdot & \cdot &  & & & \\
\cdot & & \cdot & & &  \\
\cdot &  &  & \cdot & & \\
 & & & & & 
\end{matrix}\right]
\]
and 
\[
T_2 \simeq \left[
\begin{matrix}
\begin{array}{cc|cc|}
b\, U & 0 & 0 & 0 \\
0  & 0 & 0 & 0 \\ 
\hline
0 & 0 & 0 & 0 \\
0& b_{21}I &0&0 \\ 
\hline
0&0& 0& 0 \\
0&0&0& 0  \\
0&0& b_{21^2} I &0 \\
0&0&0& b_{2^21}I  \\
\hline
\end{array} 
  &  \cdot & \cdot & \cdot & & \\

\cdot & \cdot &  & & & \\
\cdot & & \cdot & & &  \\
\cdot &  &  & \cdot & & \\
 & & & & & 
\end{matrix}\right]
\]
where $U$ is the canonical (unweighted) bilateral weighted shift operator acting on the standard
basis for $\H_\phi$. 

These decompositions can yield information on operator algebras generated by the $T_i$. For
instance, if we let $\fA_{\omega, 1}$ be the $\ca$-algebra (contained in $\B(\K_\omega)$)
generated by the $T_1$, $T_2$ from all 1-periodic shifts acting on $\K_\omega$, then we have
complete freedom on choices of scalars $a_u$, $b_v$ determining the generators of
$\fA_{\omega, 1}$. Hence, it is evident from these decompositions that $\fA_{\omega,1}$ is
unitarily equivalent to an algebra which is $\ast$-isomorphic to $\B(\H)\otimes {\rm C}(\bbT)$,
where $\H$ is a separable infinite dimensional Hilbert space and ${\rm C}(\bbT)$ is the set of
continuous functions on the unit circle (which is isomorphic to $\ca (U)$).

\vspace{0.05in}

{\noindent}${\mathbf (ii)}$ Let $T = (T_1,T_2)$ be a weighted shift on 
$\K_\omega$ of period $k=2$. We
may obtain similar block matrix decompositions for $T_1$ and $T_2$. In this case, the actions
of $T_1$, $T_2$ depend on their behaviour on the principal component $\bbF_{\omega,2}$,
which has the ordering 
\begin{eqnarray*}
 \bbF_{\omega,2}& =& \big\{ \{\phi\}, \{v\in\fnplus 1\}, \{ 2^{-1} \},  
\{v\in\fnplus 12^{-1}\}\big\} \\
 & =& \left\{ \begin{array}{ll} 
\{\phi\}, \{1\},  \{v1: v\in\fnplus, |v|=1 \}, \ldots , \\
 \{2^{-1}\}, \{12^{-1}\},  \{v12^{-1}: v\in\fnplus, |v|=1 \}, \ldots
\end{array}
\right\}.
\end{eqnarray*}
The following diagram gives the infinite tree structure here:

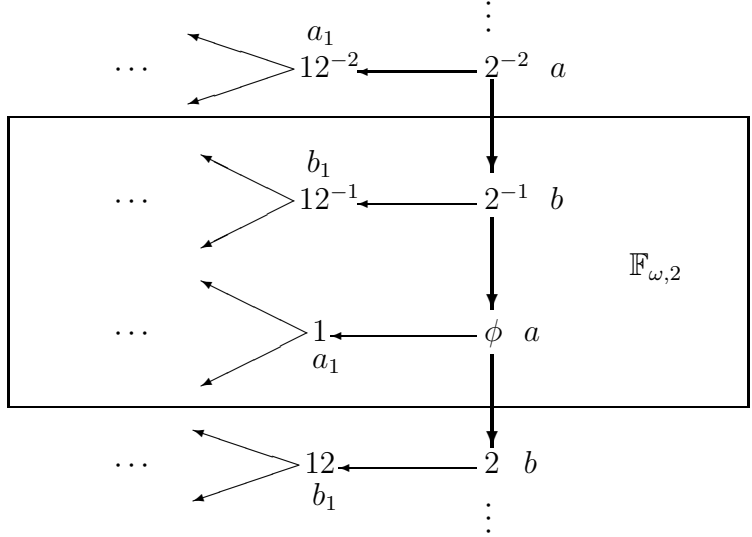
\begin{picture}(100,250)(-90,0)


\put(100,50){$2$}
\put(115,50){$b$}
\put(100,100){$\phi$}
\put(115,100){$a$}
\put(100,150){$2^{-1}$}
\put(125,150){$b$}
\put(100,200){$2^{-2}$}
\put(125,200){$a$}
\multiput(100,215)(0,5){3}{$\cdot$}
\multiput(100,35)(0,-5){3}{$\cdot$}

\put(155,125){$\bbF_{\omega,2}$}

{\thicklines      
\put(103,95){\vector(0,-1){35}}
\put(103,147){\vector(0,-1){35}}
\put(103,199){\vector(0,-1){35}}
}

\put(-80,185){\line(1,0){280}}
\put(-80,75){\line(1,0){280}}
\put(-80,75){\line(0,1){110}}
\put(200,75){\line(0,1){110}}

\put(35,100){$1$}
\put(35,90){$a_1$}
\put(97,102){\vector(-1,0){55}}

\put(32,50){$12$}
\put(35,38){$b_1$}
\put(97,52){\vector(-1,0){52}}

\put(30,150){$12^{-1}$}
\put(33,164){$b_1$}
\put(97,152){\vector(-1,0){45}}

\put(30,200){$12^{-2}$}
\put(33,214){$a_1$}
\put(97,202){\vector(-1,0){45}}

\put(33,103){\vector(-2,1){40}}
\put(33,103){\vector(-2,-1){40}}

\multiput(-40,100)(5,0){3}{$\cdot$}

\put(30,53){\vector(-3,1){40}}
\put(30,53){\vector(-3,-1){40}}
\multiput(-40,50)(5,0){3}{$\cdot$}

\put(28,153){\vector(-2,1){35}}
\put(28,153){\vector(-2,-1){35}}
\multiput(-40,150)(5,0){3}{$\cdot$}

\put(28,203){\vector(-3,1){40}}
\put(28,203){\vector(-3,-1){40}}
\multiput(-40,200)(5,0){3}{$\cdot$}

\end{picture}

From Theorem~\ref{mainthm}, we have $\K_\omega = P_0 \K_\omega \oplus P_1 \K_\omega$,
and it is easy to see that $P_0 \K_\omega$ and $P_1 \K_\omega$ are reducing subspaces for
$T_1$ in this example. Further, it is clear from the matrix decomposition associated with $T_1$,
that the restrictions $T_1|_{P_0 \K_\omega}$ and $T_1|_{P_1 \K_\omega}$ are each unitarily
equivalent to  operators in the $T_1$ class of the 1-periodic shifts 
above.  Thus
$T_1$ is unitarily equivalent to the orthogonal direct sum of two operators $T_1 \simeq
T_1^{(0)} \oplus T_1^{(1)}$ from the 1-periodic case. 
On the other hand, when we use the above ordering on $\bbF_{\omega,2}$, and make the spatial
identifications $\H_u \simeq \H_\phi$ as in the above example (where these are the subspaces
obtained in the 2-periodic case), we decompose 
\[
\K_\omega = P_0 \K_\omega \oplus P_1 \K_\omega \simeq \left( \sum_{u\in\bbG_0^{(0)}}
\oplus \H_\phi \right) \oplus  \left(
\sum_{u\in\bbG_1^{(0)}} \oplus \H_\phi \right). 
\]
For $u,v\in\bbF_{\omega,2}$, let $E_{uv}$ be the matrix unit corresponding to the $(u,v)$ entry
in this decomposition of $\B(\K_\omega)$. 
Then we may write the block matrix decomposition unitarily equivalent to $T_2$ as:
\begin{eqnarray*}
T_2 &\simeq& (b\,U) E_{2^{-1},\phi} + \sum_{u\in\bbG_0^{(0)}; \, u\neq \phi} a_{2u}
E_{2u,u} \\ 
& & + a\, E_{\phi,2^{-1}} + \sum_{u\in\bbG_1^{(0)}; \, u\neq 2^{-1}} b_{2u} E_{2u,u},
\end{eqnarray*}
with $U$ the bilateral shift acting on the standard basis for $\H_\phi$, which is 
$\{ \xi_{v_l} : v_l = 2^{2l}, l\in\bbZ\}$ in this case. 
Again, we may deduce properties of the operator algebras generated by such weighted shifts. For
instance, it follows from these matrix decompositions that the $\ca$-algebra $\fA_{\omega,2}$
is $\ast$-isomorphic to $\B(\H)\otimes {\rm C}(\bbT)$ as in the 1-periodic case. In fact, it can be
shown that all the $\ca$-algebras $\fA_{\omega,k}$, $k\geq 1$, are $\ast$-isomorphic. We
address  this point below. 
\end{eg}

\begin{eg}\label{example2}
Let $n=2$ and let $\omega = v_0v_0v_0\cdots$ belong to $\Omega_+$ with $v_0 =12$.
Consider the weighted shifts $T=(T_1,T_2)$ acting on $\K_\omega$ of period $k=2$. The
principal component, with its natural ordering, is given by 
\[
 \bbF_{\omega,2} 
  = \left\{ \begin{array}{ll} 
\{\phi\},   \{v\in\fnplus 1 \}, \{1^{-1}\},    \{v\in\fnplus 21^{-1} \}, \\
\{v_0^{-1}\},     \{v\in\fnplus 1v_0^{-1} \}, \{ (v_01)^{-1}\},    \{v\in\fnplus 2(v_01)^{-
1} \}
\end{array}
\right\}.
\]
The infinite tree structure for $\K_\omega$ traced out by a 
given 2-periodic shift $T=(T_1,T_2)$ is:

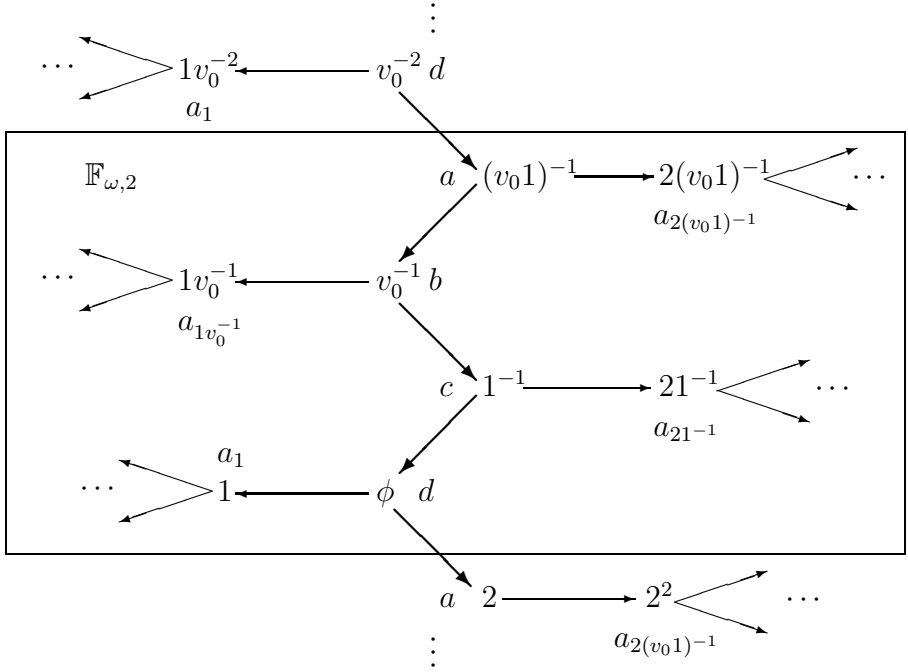
\vspace{0.3in}

\begin{picture}(260,250)(-30,0)



\multiput(130,235)(0,5){3}{$\cdot$}
\multiput(130,5)(0,-5){3}{$\cdot$}

\put(134,20){$a$}
\put(126,60){$d$}
\put(134,100){$c$}
\put(130,140){$b$}
\put(134,180){$a$}
\put(130,220){$d$}

\put(150,20){$2$}
\put(110,60){$\phi$}
\put(150,100){$1^{-1}$}
\put(110,140){$v_0^{-1}$}
\put(150,180){$(v_01)^{-1}$}
\put(110,220){$v_0^{-2}$}

{\thicklines      
\put(117,57){\vector(1,-1){29}}
\put(119,135){\vector(1,-1){29}}
\put(119,215){\vector(1,-1){29}}
\put(148,180){\vector(-1,-1){29}}
\put(148,100){\vector(-1,-1){29}}}

\put(0,180){$\bbF_{\omega,2}$}

\put(-30,40){\line(1,0){340}}
\put(-30,200){\line(1,0){340}}
\put(-30,40){\line(0,1){160}}
\put(310,40){\line(0,1){160}}

\put(35,220){$1v_0^{-2}$}
\put(38,206){$a_1$}
\put(107,223){\vector(-1,0){50}}

\put(33,224){\vector(-3,1){35}}
\put(33,224){\vector(-3,-1){35}}
\multiput(-7,222)(-5,0){3}{$\cdot$}

\put(35,140){$1v_0^{-1}$}
\put(35,125){$a_{1v_0^{-1}}$}
\put(107,143){\vector(-1,0){50}}
\put(33,144){\vector(-3,1){35}}
\put(33,144){\vector(-3,-1){35}}
\multiput(-7,142)(-5,0){3}{$\cdot$}

\put(50,60){$1$}
\put(50,75){$a_{1}$}
\put(107,63){\vector(-1,0){50}}
\put(48,64){\vector(-3,1){35}}
\put(48,64){\vector(-3,-1){35}}
\multiput(8,62)(-5,0){3}{$\cdot$}

\put(212,20){$2^2$}
\put(200,5){$a_{2(v_01)^{-1}}$}
\put(158,23){\vector(1,0){50}}
\put(223,22){\vector(3,1){35}}
\put(223,22){\vector(3,-1){35}}
\multiput(265,20)(5,0){3}{$\cdot$}

\put(217,100){$21^{-1}$}
\put(215,85){$a_{21^{-1}}$}
\put(166,103){\vector(1,0){48}}
\put(239,102){\vector(3,1){35}}
\put(239,102){\vector(3,-1){35}}
\multiput(276,100)(5,0){3}{$\cdot$}

\put(217,180){$2(v_01)^{-1}$}
\put(215,165){$a_{2(v_01)^{-1}}$}
\put(185,183){\vector(1,0){30}}
\put(257,182){\vector(3,1){35}}
\put(257,182){\vector(3,-1){35}}
\multiput(290,180)(5,0){3}{$\cdot$}



\end{picture}

\vspace{0.3in}

In this case, $\J_{\phi,2} = \{ v_0^{2l} : l\in\bbZ\}$ and $\J_{u,2} = u\J_{\phi,2}$ for
$u\in\bbF_{\omega,2}$. Thus $\K_\omega = \sum_{u\in\bbF_{\omega,2}}\oplus \H_u \simeq 
\sum_{u\in\bbF_{\omega,2}}\oplus \H_\phi$, and with the above ordering on
$\bbF_{\omega,2}$, we have the corresponding block matrix decompositions 
for $T_1$ and $T_2$: 
\[
T_1 \simeq b \,E_{v_0^{-1}, (v_01)^{-1}} + d \, E_{\phi, 1^{-1}}+
\sum_{u\in\bbF_{\omega,2}; \,\, u\notin \{1^{-1}, (v_01)^{-1}\}} a_{1u} E_{1u,u},
\]
and 
\[
T_2 \simeq (a \, U) E_{(v_0 1)^{-1},\phi}+ c \, E_{1^{-1}, v_0^{-1}} + 
\sum_{u\in\bbF_{\omega,2};\,\, u\notin \{\phi,v_0^{-1}\} } a_{2u} E_{2u,u}.
\]
where again $U$ is the canonical bilateral weighted shift operator on the standard basis for
$\H_\phi$, which is $\{ \xi_{v_l} : v_l = v_0^{2l}, l\in\bbZ\}$. 
\end{eg}

From discussions in
the previous examples, we may deduce the following result. 
Given a periodic word $\omega$ in $\Omega_+$ and a positive integer $k\geq 1$, let
$\fA_{\omega,k}$ be the $\ca$-algebra (contained in $\B(\K_\omega)$) generated by all the
$T_i$ from every $k$-periodic shift $T = \rowt$ acting on $\K_\omega$. 

\begin{thm}
Let $n\geq 2$. For every periodic word $\omega$ in $\Omega_+$ and positive integer $k\geq 1$,
$\fA_{\omega,k}$ is $\ast$-isomorphic to $\B(\H)\otimes {\rm C}(\bbT)$, where $\H$ is a
separable infinite dimensional Hilbert space. 
\end{thm}

\Prf
It is evident from the previous examples that $\fA_{\omega,k}$ is unitarily equivalent to a
$\ca$-subalgebra of $\fA = \B(\H)\otimes {\rm C}(\bbT)$. Furthermore, each matrix unit
$E_{uv}$, for $u,v\in\fomk$, is present in this subalgebra, and the general matrix
decomposition
for exactly one of the  $T_i$ has a weight multiple of the canonical bilateral shift in one of its
entries. As there is complete freedom on choices of weights for the generators of
$\fA_{\omega,k}$, it follows that this subalgebra is in fact the entire algebra $\fA$. 
\bx

\begin{rem}
This result illustrates a  difference between the commutative ($n=1$) and
noncommutative ($n\geq 2$) cases. Indeed, for $n=1$, the $\ca$-algebra generated by all 
$k$-periodic bilateral weighted shift operators, with respect to a given basis, is easily seen to be
isomorphic to the algebra $\M_k ({\rm C}(\bbT))$ of $k\times k$ matrices with entries in ${\rm
C}(\bbT)$. (These algebras played a role in the work of Bunce and Deddens 
\cite{BD1,BD2}.) Whereas 
for $n\geq 2$, the connection with $k$ is washed away by the infinite
multiplicities present, at least in this $\ca$-algebra setting. The nonselfadjoint versions of 
$\fA_{\omega,k}$ will clearly be unitarily equivalent to matrix function algebras as well;
however, we would expect distinct algebras for different values of $k$ in the nonselfadjoint
case.    
\end{rem}

\subsection{Concluding Remarks and Open Problems}\label{S:conclusion}

There are a number of open problems on the operator algebras determined by these weighted
shifts. For instance, we have not addressed the classification problem for the algebras $\fsulam$.
It would be interesting to know if $\fsulam$ can be classified by spatial means strictly in terms
of  $\omega$ and $\Lambda$, up to some sort of shift-tail equivalence classes.
\begin{prob}
Does the pair $(\omega, \Lambda)$ form a complete 
set of unitary invariants for $\fsulam$?
\end{prob}
\noindent{More} generally, we wonder about the representation theory for these algebras
\cite{Kat}, as well as reflexivity issues. 
\begin{prob}
When is $\fsulam$ reflexive, or hyper-reflexive?
\end{prob}


It would also be interesting to find a description of  the $\ca$-algebras  generated  
by these shifts, $\ca(T_1,\ldots, T_n)$. In particular, work on the 
unweighted case $\ca(S_1, \ldots, S_n)\simeq\O_n$ \cite{Renault} suggests 
the following problem. 
\begin{prob}
Can the $\ca$-algebras $\ca(T_1,\ldots,T_n)$ be described using 
groupoid techniques? 
\end{prob}

From a more operator theoretic point of view, we ask if there are 
extensions of other results on bilateral weighted shift 
operators to this setting that uncover new phenomena. 
We are also curious about the relationship between the notions of 
periodicity discovered here, and in \cite{inductive} for versions of 
unilateral weighted shifts in the noncommutative case.  
\begin{prob}
Is there a unifying framework for the noncommutative multivariable notions
of periodicity discovered here and in \cite{inductive}?   
\end{prob}


\vspace{0.1in}

{\noindent}{\it Acknowledgements.}
We would like to thank Ken Davidson for organizing a  
workshop on nonselfadjoint operator algebras at the Fields Institute 
in Toronto (July 2002), where 
the author had several illuminating discussions with participants. 
Specifically, concerning the current paper, discussions with 
Elias Katsoulis motivated us to look harder at these shifts. We are also 
grateful  to Stephen Power for helpful conversations on this topic at the workshop. Thanks to
members of the
Department of Mathematics  at Purdue University for kind hospitality during preparation of this
article.




%

\end{document}